\documentclass{article}[11pt]
\usepackage[cp1250]{inputenc}

\usepackage{authblk}
\usepackage[english]{babel}
\usepackage[colorlinks,citecolor=blue,urlcolor=blue]{hyperref}          

\usepackage{a4wide}
\usepackage{multicol}
\usepackage{amsfonts}
\usepackage{graphicx}
\usepackage{caption}
\usepackage{bm}
\usepackage{bbm} 
\usepackage{amsthm,amsmath,amssymb} 
\usepackage{color}
\usepackage{picture}
\usepackage{tikz}
\usepackage{dsfont}
\usepackage{mathrsfs}
\usepackage{textcomp}
\usepackage{mathtools}

\DeclarePairedDelimiter\floor{\lfloor}{\rfloor}
\usepackage{mleftright}
\usepackage{makecell}

\usepackage{xcolor}
\usepackage{framed}
\colorlet{shadecolor}{gray!20}

\usepackage[T1]{fontenc}

\usepackage{listings}
\usepackage{pgfplots}
\pgfplotsset{compat=1.15}

\newcommand{\R}{\mathbb{R}}
\newcommand{\N}{\mathbb{N}}

\usepackage{xargs}

\usepackage{color}
\definecolor{qqqqff}{rgb}{0,0,1}
\definecolor{ududff}{rgb}{0.30196078431372547,0.30196078431372547,1}
\definecolor{xdxdff}{rgb}{0.49019607843137253,0.49019607843137253,1}
\definecolor{iqqqqff}{rgb}{0,0,1}
\definecolor{iududff}{rgb}{1,0,0}
\definecolor{ixdxdff}{rgb}{0,1,0}

\newcommand{\set}[1]{\mathcal{#1}}
\newcommand{\one}{\mathbbm{1}}

\newcommandx{\bkq}[2]{\mathcal{B}^{#1}_{#2}}
\newcommandx{\bkqp}[2]{(\mathcal{B}^{#1}_{#2})_+}
\newcommand{\ising}{\bkqp{k}{2}}

\numberwithin{equation}{section}

\newcommand*{\lemproofname}{Proof}

\newtheorem{lemma}{Lemma}[section]
\newtheorem{theorem}[lemma]{Theorem}
\newtheorem{corollary}[lemma]{Corollary}
\newtheorem{definition}[lemma]{Definition}
\newtheorem{remark}[lemma]{Remark}
\newtheorem{example}[lemma]{Example}

\newtheorem{conjecture}[lemma]{Conjecture}
\newtheorem{proposition}[lemma]{Proposition}

\usepackage{extarrows}

\newcommand{\Then}{\Rightarrow}

\title{Level sets and maximum likelihood estimation\\ for the Ising model
}
 \author[1]{Tomasz Skalski \thanks{corresponding author, email: Tomasz.Skalski@pwr.edu.pl}}
 \author[1]{Tomasz Stroi\'nski \thanks{email: Tomasz.Stroinski@pwr.edu.pl}}
 \affil[1]{Faculty of Pure and Applied Mathematics, Wroc\l{}aw University of Science and Technology, Wybrze\.{z}e Wyspia\'{n}skiego 27, 50-370 Wroc\l{}aw, Poland}
 \date{\today}

\providecommand{\keywords}[3]{\small\textbf{Key words:} #1 \\
\indent \textbf{Mathematics Subject Classification (#2):} #3}

\begin{document}
\maketitle

\begin{abstract}
\noindent
Bogdan et al. established a new criterion to determine the~existence of a maximum likelihood estimator in discrete exponential families. It uses the notion of the set of uniqueness, which allows one to apply the~problem to the Ising model from statistical mechanics. We propose a~full~characterization of the existence of the MLE for the Ising model among the~level~sets used in related combinatorial problems. We then establish new bounds for the size of the smallest set of uniqueness for the products of~Rademacher functions.
\end{abstract}

\keywords{maximum likelihood, discrete exponential family, Ising model, level sets.}{2010}{62F10, 05D05.}

\section{Introduction}

The class of exponential families is one of the most popular classes of distributions in probability and statistics. Apart from many properties which make them useful to derive better theoretical results and apply in the practical setting~\cite{MR754971}, \cite[Sec. 2.7]{MR1639875}, they facilitate the analysis of the Maximum Likelihood Estimator (MLE). However, in some cases the MLE may not exist, as pointed out, e.g., by Crain~\cite{MR0362678, MR0428559} in situations where the number of parameters is too large. Barndorff-Nielsen~\cite[Theorem 9.13]{MR489333} established a full characterization of~the~existence of the MLE for exponential families, however, in general it is difficult to apply in practice~\cite{MLE}. Bogdan et al.~\cite{MLE} proposed a new, straightforward criterion for the existence of the MLE for discrete exponential families, which was inspired by the article by~Bogdan and Bogdan~\cite{MR1768245} in the continuous case. This criterion uses the notion of~the~set~of~uniqueness, i.e., the set which ensures that the function from a given class which vanishes on such a set is the constant zero function on the whole space.

The article by Bogdan et al.~\cite{MLE} enriches the discussion on the existence of the~MLE with more exact results on the Rademacher functions and their products, both of~which are related to the discrete hypercube $\set{X}=\{-1,+1\}^k$. These results include the theoretical conditions for the existence of~the~MLE for the given sample, as~well as~the~probability bounds for the sample size needed to ensure it. However, in general, for the products of Rademacher functions, establishing the~exact conditions ensuring the existence of the MLE remains an open problem.

In this paper, we focus on the fact that in many cases of the products of Rademacher functions, the minimal samples for which the MLE exists, can be described as the~level sets, i.e., the sets of those points in $\set{X}=\{-1,+1\}^k$, for which the number of~positive coordinates belongs to a given set of numbers. The level sets have been used, e.g., for establishing the smallest subsets of $\set{X}$ which intersect every subcube of the given size~\cite{MR2299695},~\cite[Sec. 2.5]{baber},~\cite{MR1204793}. This problem is also known as the subcube allocation problem~\cite{MR1204793}. It can be used, e.g., in the analysis of the hypercube processors, known from their interconnection properties, their fault-tolerance as well as the ability to support multiple users, each assigned a subcube of~a~given size, cf.~\cite{DuttHayes}.

The paper is composed as follows. Section~\ref{sec:level} includes the full characterization of the existence of the MLE for the level sets in the Ising model. In Section~\ref{sec:extr} we~discuss the bounds on the size of the smallest set of uniqueness for the products of~Rademacher functions, improving the upper bounds. Section~\ref{sec:proofs} covers the proofs from Sections~\ref{sec:ising} and~\ref{sec:level}. 

\subsection{Discrete exponential families}
For $k \in \N$, let us consider $\set{X}:=\{-1, +1\}^{k}$, the $k$-dimensional discrete cube with the uniform weight
$\mu(x) = 2^{-k}$, $x \in \set{X}$. Thus, $K=|\set{X}|=2^k$. We adopt the convention of~\cite{MLE} that $\R^{\set{X}}$ is the class of all real-valued functions on $\set{X}$. Throughout the paper, we let $\set{B}$ denote a fixed linear subspace of $\R^\set{X}$ containing the constant function $\varphi\equiv 1$.

We often refer to the cone of non-negative functions $\set{B}_+ := \{\phi \in \set{B}: \phi \geq 0\}$.\\
We define, as in~\cite{MLE}, the exponential family
$e(\set{B}) = \{e(r): r \in \set{B}\}$, where\\ $e(r)=\exp\{r-\psi(r)\}$ is an exponential density and 
$$
\psi(r)=\log\left(\sum\limits_{x\in\set{X}}e^{r(x)}\mu(x)\right)
$$
is the log-partition function.

As $\set{X}$ is a finite set, $e(\set{B})$ is known as a discrete exponential family, cf.~\cite{MLE}.\\
Let $x_1,\ldots,x_n\in\set{X}$ be an iid sample. For an exponential density $p=e(r)$, its likelihood function $L_p(x_1,\ldots,x_n)$ 
is defined as $p(x_1)\cdot \ldots\cdot p(x_n)$ 
and its log-likelihood function $l_p(x_1,\ldots,x_n)$ is defined as $\log L_p(x_1,\ldots,x_n)$. \\
An element $\hat{p}$ of $e(\set{B})$ is called a maximum likelihood estimator (MLE)\\ for~$x_1,\ldots,x_n$ and~$e(\set{B})$ if $$
L_{\hat{p}}\left(x_1,\ldots, x_n\right) = \sup_{p \in e(\set{B})} L_{p}\left(x_1,\ldots, x_n\right),$$ 
or, equivalently,
$$
l_{\hat{p}}\left(x_1,\ldots, x_n\right) = \sup_{p \in e(\set{B})} l_{p}\left(x_1,\ldots, x_n\right).
$$
To establish new elegant results on the existence of the MLE for the exponential family $e(\set{B})$, we use the notion of the set of uniqueness from~\cite{MLE}. 
\begin{definition}(\cite[\rm{Sec. 2.}]{MLE})
Let $U\subset \set{X}$. Let $\set{B}$ be a fixed linear subspace of $\R^{\set{X}}$. 
We say that $U$ is a set of uniqueness for $\set{B}$ if $\phi=0$ is the only function in~$\set{B}$ such that $\phi = 0$ on $U$. Similarly, we say that $U$ is a set of uniqueness for $\set{B}_+$ if~$\phi=0$ is the only function in $\set{B}_+$ such that $\phi = 0$ on $U$.
\end{definition}
\begin{remark}(\cite[\rm{Sec. 2.}]{MLE})
If $U$ is a set of uniqueness for $\set{B}$, then $U$ is a set of~uniqueness for $\set{B}_+$, since $\set{B}_+ \subset \set{B}$.
\end{remark}

The main motivation behind the sets of uniqueness is the following criterion\\ for the existence of the MLE, stated by Bogdan et al.~\cite{MLE}.
\begin{theorem}(\cite[\rm{Theorem 2.2.}]{MLE})\label{thm:11}
MLE for $e(\set{B})$ and $x_1,\ldots,x_n\in\set{X}$ exists if and only if $\{x_1,\ldots,x_n\}$ is a set of uniqueness for $\set{B}_+$.
\end{theorem}
\subsection{Rademacher functions}
\label{sec:rad}
Theorem~\ref{thm:11} motivated us to study sets of uniqueness for the space of linear combinations of Rademacher functions and products of pairs of Rademacher functions. Let us recall that Rademacher functions\\ $r_j\ :\ \set{X} \to \R$, $j=0,\ldots,k$, is defined by $r_0\equiv 1$, and
\begin{align*}
r_j(x)=x_j, && {\rm{for\ }} x=(x_1,\ldots,x_k)\in\set{X}
\end{align*}
for $j=1,\ldots,k$. Let $\set{B}^k = \mbox{Lin}\{r_0, r_1, \ldots, r_k\}$.\\
A $(k-q)$-subcube of $\set{X}$ is defined as a subset of $\set{X}$ formed by~fixing exactly $q$ coordinates at $+1$ or $-1$, that is for a given
$$
{\bf{\varepsilon}}=\left( \varepsilon_1,\varepsilon_2,\ldots,\varepsilon_k\right)\in\{-1,+1,*\}^{k}
$$
such that $|\{\varepsilon_j = *\}| = k-q$, the corresponding $(k-q)$-subcube is 
\begin{align*}
\{x\in\set{X}:\ x_j = \varepsilon_j, && {\rm{for\ all}}&& j\in\{1,\ldots,k\} && {\rm{such\ that}}&& \varepsilon_j \neq *\}.
\end{align*}
Here $*$ denotes a coordinate which is not fixed.\\ 
Observe that each $(k-q)$-subcube consists of exactly $2^{k-q}$ points.  For~example, for $k=4$, the $2$-subcube defined by~fixing the first coordinate at $+1$ and~the~last coordinate at $-1$ is denoted as~$(+1,*,*,-1)$.
With the notion of the $(k-q)$-subcube we discuss the spaces spanned by products of $q$ Rademacher functions.
Let $k \in \N$ and $1 \leq q \leq k$. We define
\begin{equation*}
\set{B}_q^k =  
{\mbox{Lin}\left\{w_L: L \subset \left\{1,\ldots, k\right\} \mbox{ and } |L| \leq q\right\},}
\end{equation*}
where 
$$
w_{L}({x}) = \prod_{j\in{L}}{r_j}({x}),\quad {x}\in{\{-1,+1\}^k},\quad L\subset\{1,\ldots,k\},
$$
are the Walsh functions, see, e.g., Bogdan et al.~\cite{MLE} or Jendrej et al. \cite{MR3446023}.
Similarly, we define the non-negative cone of functions in $\bkq{k}{q}$ as $\bkqp{k}{q}$.


\begin{remark}(\cite[\rm{Remark~5.6}]{MLE})
\label{rmk:SoU_Bkq}
Let $1\leq q_1 \leq q_2 \leq k$. Then every set $U$ of~uniqueness for $(\set{B}_{q_2}^k)_+$ is of uniqueness for 
$(\set{B}_{q_1}^k)_+$, since 
$(\set{B}_{q_1}^k)_+ \subset (\set{B}_{q_2}^k)_+$.\\
Analogously, every set $U$ of uniqueness for $\set{B}_{q_2}^k$ is of uniqueness for 
$\set{B}_{q_1}^k$, since 
$\set{B}_{q_1}^k \subset \set{B}_{q_2}^k$.
\end{remark}
In this paper, we primarily focus on the space spanned by the products of two Rademacher functions, $\bkq{k}{2}$, and its non-negative cone, $\ising$, which corresponds to the~Ising model described below.
\subsection{Ising model}~\label{sec:ising}
The Ising model was proposed by Lenz in 1920 and developed by Ising in 1925~\cite{brush}. 
\begin{definition}(\cite{MR2650042})
\label{def:boltzmann}
Let $G=(V,E)$ be an undirected graph, where\\ $V\!=\!\{1,2,\ldots,k\}$ is a set of its~vertices and $E\subset\binom{V}{2}=\{(i,j):1\leq\! i\! <\! j\! \leq\! k\}$ is~a~set of its edges. Let ${{x}} = \{x_i:\ i\in V\}$, where $x_i \in \{-1, +1\}$. 
The Ising model on $G$ is defined as the probability space $(\set{X}, \mu_{B,\beta})$, where $\set{X}=\{-1,+1\}^V$ is the sample space and $\mu_{B,\beta}$ is the Boltzmann distribution:
\begin{equation*}
\mu_{B,\beta}({{x}}) = \frac{1}{Z(B,\beta)} \exp\left(B\ \sum\limits_{i\in V} x_i + \beta\ \sum\limits_{(i,j)\in E} x_i x_j \right),
\end{equation*}
where 
$$Z(B,\beta)\ = \sum\limits_{{x}\in\{-1,+1\}^V}{\exp\left(B\sum\limits_{i\in V} x_i + \beta\sum\limits_{(i,j)\in E}x_i x_j \right) }$$ 
is the normalizing constant (partition function). The parameters $B$ and $\beta$\\ can be interpreted, respectively, as the magnetic field and the inverse temperature.
\end{definition}

In our case, we assume the graph $G=(V,E)$ to be a complete graph consisting of $k$ vertices. Thus, each ${x}$ is a point in $\{-1,+1\}^k$, every vertex corresponds to a single Rademacher function and each pair of vertices (i.e. each edge) corresponds to a product of two different Rademacher functions. The constant function $f\equiv 1$ is implicitly hidden by the partition function.\\ 
When a graph $G=(V,E)$ is not complete, the function space $\set{B}$ is spanned by the constant function $f\equiv 1$, all single Rademacher functions and the products of pairs of Rademacher functions which correspond to the pairs $(i,j)$ for which $(i,j)\in E$.

\begin{example}
\label{ex:ising}
\leavevmode
\makeatletter
\@nobreaktrue
\makeatother
\begin{enumerate}
    \item If $G$ is an empty graph, i.e., $(V,E)=(\{1,2,3\},\emptyset)$, then $\set{B}_+=(\set{B}_1^3)_+$. 
    \item Let $G'$ consist of three vertices and one edge, i.e., $V'=\{1,2,3\}$ and ${E'=\{(1,2)\}}$. Then $\set{B}'$ is spanned by $f\equiv 1$, all single Rademacher functions and the function $r_1 r_2$. In terms of subcubes, it is the linear span of the characteristic functions of all faces and the edges parallel to the edge $(-1,-1,*).$
    \item Let $G''$ consist of three vertices and two edges, i.e., $V''=\{1,2,3\}$ and $E''=\{(1,2),(1,3)\}$. Then $\set{B}''$ is spanned by $f\equiv 1$, all single Rademacher functions and the functions $r_1 r_2$ and $r_1 r_3$. In terms of subcubes, it is the linear span of the characteristic functions of all faces and the edges parallel to the edge $(-1,-1,*)$ or to $(-1,*,-1)$. 
    \item If $G$ is a complete graph, i.e., $E=\binom{V}{2}$, then $\set{B}_+=\bkqp{3}{2}$.
\end{enumerate}
\end{example}
\begin{figure}
    \centering
    \begin{tikzpicture}[line cap=round,line join=round,x=1cm,y=1cm]
\draw [line width=0.8pt,color=iqqqqff] (3,2)-- (3,6);
\draw [line width=0.8pt,color=iududff,dotted] (3,2)-- (7,2);
\draw [line width=2.8pt,color=iqqqqff] (7,2)-- (7,6);
\draw [line width=2.8pt,color=iududff, dash pattern=on 1pt off 4pt] (3,6)-- (7,6);
\draw [line width=0.8pt,color=iqqqqff] (0,0)-- (3,2);
\draw [line width=2.8pt,color=iqqqqff] (4,0)-- (7,2);
\draw [line width=2.8pt,color=iqqqqff] (0,4)-- (3,6);
\draw [line width=2.8pt,color=iqqqqff] (4,4)-- (7,6);
\draw [line width=2.8pt,color=iqqqqff] (0,0)-- (0,4);
\draw [line width=2.8pt,color=iududff, dash pattern=on 1pt off 4pt] (0,0)-- (4,0);
\draw [line width=2.8pt,color=iqqqqff] (4,0)-- (4,4);
\draw [line width=2.8pt,color=iududff, dash pattern=on 1pt off 4pt] (0,4)-- (4,4);

\draw [color=qqqqff](-0.6,-0.1) node[anchor=north west] {$(-1,-1,-1)$};
\draw [color=qqqqff](3.4,-0.1) node[anchor=north west] {$(-1,-1,+1)$};
\draw [color=qqqqff](-0.1,3.9) node[anchor=north west] {$(-1,+1,-1)$};
\draw [color=qqqqff](3.9,3.9) node[anchor=north west] {$(-1,+1,+1)$};
\draw [color=qqqqff](2.4,1.9) node[anchor=north west] {$(+1,-1,-1)$};
\draw [color=qqqqff](6.4,1.9) node[anchor=north west] {$(+1,-1,+1)$};
\draw [color=qqqqff](2.9,5.9) node[anchor=north west] {$(+1,+1,-1)$};
\draw [color=qqqqff](6.9,5.9) node[anchor=north west] {$(+1,+1,+1)$};
\begin{scriptsize}
\draw [fill=qqqqff] (4,0) circle (2.5pt);
\draw [fill=qqqqff] (0,0) circle (2.5pt);
\draw [fill=qqqqff] (0,4) circle (2.5pt);
\draw [fill=qqqqff] (4,4) circle (2.5pt);
\draw [fill=qqqqff] (7,2) circle (2.5pt);
\draw [fill=qqqqff] (3,2) circle (2.5pt);
\draw [fill=qqqqff] (3,6) circle (2.5pt);
\draw [fill=qqqqff] (7,6) circle (2.5pt);
\end{scriptsize}
\end{tikzpicture}
\caption{Illustration of Example~\ref{ex:ising}.2. The characteristic functions of the red (horizontal, dotted) edges are in $\set{B'}$.}
\label{fige2}
\end{figure}

\begin{figure}
    \centering
    \begin{tikzpicture}[line cap=round,line join=round,x=1cm,y=1cm]
    
\draw [line width=0.8pt,color=ixdxdff, dashed] (3,2)-- (3,6);
\draw [line width=0.8pt,color=iududff, dotted] (3,2)-- (7,2);
\draw [line width=2.8pt,color=ixdxdff, dash pattern=on 12pt off 4pt] (7,2)-- (7,6);
\draw [line width=2.8pt,color=iududff, dash pattern=on 1pt off 4pt] (3,6)-- (7,6);
\draw [line width=0.8pt,color=iqqqqff] (0,0)-- (3,2);
\draw [line width=2.8pt,color=iqqqqff] (4,0)-- (7,2);
\draw [line width=2.8pt,color=iqqqqff] (0,4)-- (3,6);
\draw [line width=2.8pt,color=iqqqqff] (4,4)-- (7,6);
\draw [line width=2.8pt,color=ixdxdff, dash pattern=on 12pt off 4pt] (0,0)-- (0,4);
\draw [line width=2.8pt,color=iududff, dash pattern=on 1pt off 4pt] (0,0)-- (4,0);
\draw [line width=2.8pt,color=ixdxdff, dash pattern=on 12pt off 4pt] (4,0)-- (4,4);
\draw [line width=2.8pt,color=iududff, dash pattern=on 1pt off 4pt] (0,4)-- (4,4);

\draw [color=qqqqff](-0.6,-0.1) node[anchor=north west] {$(-1,-1,-1)$};
\draw [color=qqqqff](3.4,-0.1) node[anchor=north west] {$(-1,-1,+1)$};
\draw [color=qqqqff](-0.1,3.9) node[anchor=north west] {$(-1,+1,-1)$};
\draw [color=qqqqff](3.9,3.9) node[anchor=north west] {$(-1,+1,+1)$};
\draw [color=qqqqff](2.4,1.9) node[anchor=north west] {$(+1,-1,-1)$};
\draw [color=qqqqff](6.4,1.9) node[anchor=north west] {$(+1,-1,+1)$};
\draw [color=qqqqff](2.9,5.9) node[anchor=north west] {$(+1,+1,-1)$};
\draw [color=qqqqff](6.9,5.9) node[anchor=north west] {$(+1,+1,+1)$};
\begin{scriptsize}
\draw [fill=qqqqff] (4,0) circle (2.5pt);
\draw [fill=qqqqff] (0,0) circle (2.5pt);
\draw [fill=qqqqff] (0,4) circle (2.5pt);
\draw [fill=qqqqff] (4,4) circle (2.5pt);
\draw [fill=qqqqff] (7,2) circle (2.5pt);
\draw [fill=qqqqff] (3,2) circle (2.5pt);
\draw [fill=qqqqff] (3,6) circle (2.5pt);
\draw [fill=qqqqff] (7,6) circle (2.5pt);
\end{scriptsize}
\end{tikzpicture}
\caption{Illustration of Example~\ref{ex:ising}.3. The characteristic functions of the red (horizontal, dotted) and green (vertical, dashed) edges are in $\set{B''}$.}
\label{fige3}
\end{figure}

\begin{remark}~\label{rem:graphs}
Let $G=(V,E)$ and $\tilde{G}=(V,\tilde{E})$, where $E\subset\tilde{E}$. Let $\set{B}$ and $\tilde{\set{B}}$ be the function spaces of the Ising model on $G$ and $\tilde{G}$, respectively. Then, by Remark~\ref{rmk:SoU_Bkq}, every set of uniqueness for $\tilde{\set{B}}_+$ is a set of uniqueness for $\set{B}_+$, since $\set{B}_+ \subset \tilde{\set{B}}_+$.
\end{remark}

To convince the reader that the Ising models from Example~\ref{ex:ising} are not the same, we state the following remark.
\begin{remark}\label{rem:compare}
Let $k=3$. Consider the Ising models from Example~\ref{ex:ising}.
\begin{enumerate}
    \item The set $\{(-1,-1,-1),(+1,+1,+1)\}$ is of uniqueness for $\bkqp{3}{1}$, but not for $\set{B}'_+$.
    \item The set $\{(+1,-1,-1), (-1,+1,-1), (-1,-1,+1), (-1,+1,+1)\}$ is of uniqueness for $\set{B}'_+$, but not for $\set{B}''_+$.
    \item The set $\set{X}\setminus\{(-1,-1,-1),(+1,+1,+1)\}$ is of uniqueness for $\set{B}''_+$, but not for $\bkqp{3}{2}$.
\end{enumerate}
\end{remark}
The proofs are given in Section~\ref{sec:proofs}.

The discrete exponential family approach to the Ising model was applied, e.g., by Johnson et al.~\cite[Sec. 2.2]{MR3595149}, Wainwright and Jordan~\cite[Example 3.1]{Ising} and~Wang et al.~\cite{MR3911110}. Wang et al.~\cite{MR3911110} treat the Ising model as a subclass of discrete hierarchical models obtained by choosing a set of rows from the Walsh matrix, where each row corresponds to the product of Rademacher functions $r_i$ for every $i$ from $I\subset\{1,\ldots,k\}$. For the Ising model, they choose the rows related to~the~products of at most two Rademacher functions, which generate exactly the $\ising$ space, cf.~\cite[Sec. 5]{MLE}.

One may observe that the partition function in the Ising model is the sum of $2^k$ elements, which makes it intractable for larger values of $k$~\cite{MR3899777}. However, even in this case the existence of the MLE is still relevant because of various methods being used to approximate the likelihood, including computation of the pseudolikelihood~\cite{besag,hofling} or the Monte Carlo methods~\cite{MR3899777}.

\section{Level sets}
\label{sec:level}
\subsection{Introduction and motivation}
Since $|\set{X}|=|\{-1,+1\}^k|=2^k$, there are exactly $2^{2^k}$ subsets of $\set{X}$. Thus, establishing all sets~of~uniqueness having the smallest cardinality becomes a challenging problem.\\
Therefore, we consider a specific class of~subsets of~$\set{X}=\{-1,+1\}^k$, called~level sets. For a~given ${x}\in\set{X}$ and $D\subset\{0,1,\ldots,k\}$, the level set $W_D^k({x})$ is defined as~the~set of all points which differ from ${x}$ at~exactly $d$~coordinates for at least one $d\in D$. To~be~more specific, for ${x}$ and $D$ we define the level set as:
\begin{equation}
\label{eq:levelset}
W_D^k({x}) = \left\{{y}\in\{-1,+1\}^k:\ {\rm{dist}}({x},{y})\in D\right\}, 
\end{equation}
where $\rm{dist}$ is the Hamming distance, i.e., 
$$
{\rm{dist}}(x,y):=|\{i:\ x_i \neq y_i\}|.
$$

Without loss of generality, we denote $W_D^k$ by $W_D^k((-1,\ldots,-1))$, which is isomorphic to $W_D^k(x)$ for every $x\in\set{X}$. When $k$ is known from the context, we simplify the notation to $W_D$. When $D$ contains only one point, i.e., $D=\{i\}$, we denote $W_{\{i\}}$ as $W_i$.

The optimality of level sets appears in a merely related problem, that is a coloring of $\set{X}=\{-1,+1\}^k$ with $k-q$ colors such that every $(k-q)$-subcube claims all $k-q$ colors. As proven in~\cite[Theorem 7]{MR2299695}, the~asymptotically maximal number of colors $c$, for which such a coloring is admissible, equals $k-q+1$. The proof is algorithmic and shows that the optimal distribution of~$c$ colors onto the vertices of $\set{X}$ is to give the color $j=0,1,\ldots,k-q$ to~the~level~set $W_{\{j,j+c,j+2c,\ldots\}}$.\\ 
In the same paper, the authors connect the $(k-q)$-chromatic coloring with the smallest cardinality of the subset of $\set{X}$ intersecting its every $(k-q)$-subcube. In fact, they denote such cardinality as $g(k,k-q)$ and prove it to be not greater than ${2^k}/(k-q+1)$.\\
Level sets often appear as important elements of minimal sets of uniqueness for $(\set{B}^k_q)_+$ with respect to inclusion, as~summarized below. 
\begin{remark}
\label{rem:012}
\leavevmode
\makeatletter
\@nobreaktrue
\makeatother
\begin{enumerate}
    \item Let $q=k$. Then the only set of uniqueness for $(\set{B}^k_k)_+$ is the full space $W_{\{0,1,2,\ldots,k\}} = \set{X}$.
    \item Let $q=k-1$. Then the only two minimal (by inclusion) sets of uniqueness for $(\set{B}^k_{k-1})_+$ are isomorphic, i.e., $W_{\{0,2,4,\ldots\}} \cong W_{\{1,3,5,\ldots\}}$.
    \item Let $q=1$. Then the smallest set of uniqueness for $(\set{B}^k_1)_+$ is the set $W_{\{0,k\}}$.
\end{enumerate}
\end{remark}
\begin{proof}
The first statement was proven in~\cite[Lemma 3.1]{MLE}. The second statement follows directly from~\cite[Theorem 5.3]{MLE}. 
Regarding the third statement, $W_{\{0,k\}}$ consists of exactly two points $x=(-1,-1,\ldots,-1)$ and~$-x=(1,1,\ldots,1)$. It is~a~set of uniqueness for $(\set{B}_1^k)_+$, which was proven in~\cite[Example 3.7]{MLE}. Smaller sets consist only of one or zero points, therefore none of~them intersects every half-cube of $\set{X}$. This contradicts~\cite[Theorem 3.6]{MLE}, and~therefore the sets isomorphic to $W_{\{0,k\}}$ are the only sets of uniqueness for $(\set{B}_1^k)_+$ having the smallest cardinality possible.
\end{proof}

To derive more complex results on the sets of uniqueness, we need the following representation. Let $\set{S}^k_q$ be the set of all $(k-q)$-subcubes of $\set{X}$. That said, each function from $\bkq{k}{q}$ can be denoted as 
\begin{align*}
\varphi(x) = \sum\limits_{S\in\set{S}^k_q}\alpha_S \one_{\{x\in S\}},
\end{align*}
where $\alpha_S\in\R$. 
Recall that every $(k-q)$-subcube is formed by fixing exactly $q$ coordinates either at $+1$ or at $-1$. Notice that the number of fixed positive coordinates may vary from $0$ to $q$. For every $0\leq j\leq q$, let $S_j$ be the set of $(k-q)$-subcubes with exactly $j$ fixed positive coordinates.
Therefore, we~perform the partition 
\begin{equation}\label{eq:partition}
\set{S}^k_q = S_0 \cup S_1 \cup \ldots \cup S_q.
\end{equation}

\subsection{Full characteristics for the Ising model}
We start with a useful lemma, which states that the sum of values of a function $\varphi\in\bkq{k}{q}$ over $\set{X}$ 
can~be~represented by $q+1$ coordinates $T_0,T_1,\ldots,T_q$, where 
\begin{equation}\label{eq:Tdef}
T_j = \sum_{S\in S_j} \alpha_S\qquad \mathrm{for\ every}\qquad j=0,1\,\ldots,q
\end{equation}
and $\alpha_S\in\R$ for every $S\in\set{S}^k_q = S_0 \cup S_1 \cup \ldots \cup S_q.$ 
\begin{lemma}\label{lem:weights:q}
Let $\varphi\in\set{B}^k_q$. Then
\begin{equation}\label{eq:T}
\sum\limits_{x\in \set{X}}\varphi(x) = 2^{k-q} \sum\limits_{j=0}^{q}T_j.
\end{equation}
\end{lemma}
\begin{proof}
\begin{align*}
&\sum\limits_{x\in\set{X}}\varphi(x) = \sum\limits_{x\in\set{X}}\sum\limits_{S\in \set{S}^k_q} \alpha_S \one_{\{x\in S\}} = \sum\limits_{S\in \set{S}^k_q} \alpha_S \sum\limits_{x\in\set{X}} \one_{\{x\in S\}}.
\end{align*}
The partition~\eqref{eq:partition} allows splitting the sum over $\set{S}^k_q$ into the sums $S_0, S_1, \ldots, S_q$. Therefore
\begin{align*}
&\sum\limits_{x\in\set{X}}\varphi(x) = \sum\limits_{S\in S_0}\alpha_S \sum\limits_{x\in\set{X}}\one_{\{x\in S\}}\! +\! \ldots\! +\! \sum\limits_{S\in S_q}\alpha_S \sum\limits_{x\in\set{X}}\one_{\{x\in S\}}\\
&= \sum\limits_{S\in S_0}\alpha_S |S|\! +\! \ldots\! +\! \sum\limits_{S\in S_q}\alpha_S |S|
\end{align*}
Every subcube $S$ has the same size $2^{k-q}$, which implies that
\begin{align*}
&\sum\limits_{x\in\set{X}}\varphi(x)
= \sum\limits_{S\in S_0}\alpha_S 2^{k-q}\! +\! \ldots\! +\! \sum\limits_{S\in S_q}\alpha_S 2^{k-q}\\
&= 2^{k-q}\left(\sum\limits_{S\in S_0} \alpha_S + \ldots + \sum\limits_{S\in S_q} \alpha_S \right)\\ &= 2^{k-q}\left( T_0 + T_1 + \ldots + T_q \right).
\end{align*}
\end{proof}
Therefore, showing that $T_0 = \ldots = T_q = 0$ is sufficient to prove that $\varphi\equiv 0$.\\
Again, we consider the partition~\eqref{eq:partition} of the set $\set{S}^k_q$ of the $(k-q)$-subcubes.
\begin{lemma}\label{lem:card:weights}
Let $S\in S_i\subset\set{S}^k_q$. Then $|W_j \cap S| = \binom{k-q}{j-i}$, with~the~convention
$\binom{n}{r}=0$ for $r<0$ and for $r>n$.
\end{lemma}
\begin{proof}
By definition, every point from $S\in S_i$ has $q$ fixed coordinates, where $i$ of them are positive and~$q-i$ are negative. To obtain a point from $S$ having exactly $j$ positive coordinates in total, one needs to~choose exactly $j-i$ from $k-q$ remaining ones, which can be done in $\binom{k-q}{j-i}$ ways.
\end{proof}

\begin{lemma}\label{rem:three}
Let $k\geq3$ and $q=2$. Then the set $W_{\{1,k\}}$ is a minimal set of uniqueness for $(\set{B}^k_2)_+$ with respect to inclusion.
\end{lemma}
\begin{proof}
Let $\varphi\in(\set{B}^k_2)_+$. 
Then
$\varphi$ is non-negative and can be denoted as
\begin{align*}
\varphi(x) &= \sum\limits_{S\in\set{S}^k_2} \alpha_S \one_{\{x\in S\}}\\
&= \sum\limits_{S\in S_0} \alpha_S \one_{\{x\in S\}} + \sum\limits_{S\in S_1} \alpha_S \one_{\{x\in S\}} + \sum\limits_{S\in S_2} \alpha_S \one_{\{x\in S\}},
\end{align*}
where $\alpha_S\in\R$ and $x\in\set{X}$. As $\varphi$ is non-negative, to show that $\varphi\equiv 0$ on $\set{X}$ it~is~sufficient to show that 
the sum of $\varphi(x)$ over $\set{X}$ equals $0$.
By Lemma~\ref{lem:weights:q}, this is equivalent to ${T_0 + T_1 + T_2 = 0}$.\\
Assume that $\varphi(W_{\{1,k\}})=\{0\}$. Since $W_1\subset W_{\{1,k\}}$, 
the sum of $\varphi(x)$ over $W_1$ equals $0$. On the other hand,
\begin{align*}
\sum\limits_{x\in W_1}\varphi(x) &=  \sum\limits_{S\in S_0} \alpha_S |W_1 \cap S| + \ldots +  
\sum\limits_{S\in S_2} \alpha_S |W_1 \cap S|\\ 
&= \binom{k-2}{1-0} \sum\limits_{S\in S_0} \alpha_S + \binom{k-2}{1-1} \sum\limits_{S\in S_1} \alpha_S + \binom{k-2}{1-2} \sum\limits_{S\in S_2} \alpha_S\\ &= (k-2) T_0 + T_1.
\end{align*}
Moreover, $(+1,+1,\ldots,+1)\in W_{\{1,k\}}$, thus $\varphi\left(\left(+1,+1,\ldots,+1\right)\right) = 0$.\\ On the other hand, 
 $(+1,+1,\ldots,+1)\in S_2$. Therefore
\begin{equation*}
\varphi\left(\left(+1,+1,\ldots,+1\right)\right) = \sum\limits_{S\in S_2} \alpha_S = T_2.
\end{equation*}
It implies that $T_2 = 0$. Now we want to show that $T_0$ and $T_1$ are non-negative.\\ Using the non-negativity of the function $\varphi$, we obtain
\begin{equation*}
0 \leq \varphi\left(\left(-1,-1,\ldots,-1\right)\right) = \sum\limits_{S\in S_0} \alpha_S = T_0,
\end{equation*}
and
\begin{align*}
0 \leq \sum\limits_{x\in W_{k-1}}\varphi(x) &= 0 + \sum\limits_{S\in S_1} \alpha_S |W_{k-1} \cap S| + 
\sum\limits_{S\in S_2} \alpha_S |W_{k-1} \cap S|\\ &= T_1 + (k-2) T_2.
\end{align*}
Since $T_2$ equals zero, we also have $T_1 \geq 0$.
Thus, considering the previous calculations, $T_0 = T_1 = T_2 = 0$. We deduce that $\varphi(\set{X})=\{0\}$ and $W_{\{1,k\}}$ is a set of~uniqueness for~$\left(\set{B}^k_2\right)_+$.\\
Now it suffices to observe that no smaller set, with respect to inclusion, is of uniqueness for $\ising$. Indeed, if $U\subset W_{\{1,k\}}$, then we have one of the following situations.\\
If $U\subset W_1$, then there exists the function $\varphi(x)=\one_{\{x\in(+1,+1,*,*,\ldots,*)\}}\in\ising$, which vanishes on~$U$, since no point with its first two coordinates positive has only one positive coordinate.\\
On the other hand, if $U\not\subset W_1$, then there exists $1\leq j\leq k$ such that $$(-1,\ldots,-1,\underbrace{+1}_{j},-1,\ldots,-1)\notin U.$$
Then the function $\varphi(x)=\one_{\{x\in S\}}\in\ising$ with $S\in S_1$ having its $j^{\rm{th}}$ coordinate fixed as $+1$, vanishes on $U$. \\ Therefore $W_{\{1,k\}}$ is a minimal set of~uniqueness for~$\left(\set{B}^k_2\right)_+$.
\end{proof}
Our main result is the full characterization of sets of uniqueness for $\ising$ among the level sets in $\set{X}=\{-1,+1\}^k$.

\begin{theorem}\label{thm:main}
Let $k\geq 3$ and $D\subset\{0,1,\ldots,k\}$. The set $W_{D}$ is of uniqueness for~$(\set{B}^k_2)_+$ if and only if there exist $i,j\in D$ such that $2\leq |i-j|\leq k-1$.
\end{theorem}

The proof of this theorem is geometric. It uses the partition~\eqref{eq:partition} of the set of all $(k-2)$-subcubes of $\set{X}$ as well as the mapping of all level sets $W_j$: $j=0,1,\ldots,k$ onto the corresponding points $P^{(j)}\in\R^2$ forming a convex polytope. Later on, the~convexity of this polytope is applied. The proof is given in Section~\ref{sec:proof:main}.

By Remark~\ref{rem:graphs} one may observe that any set $W_D$ which is of uniqueness for $\ising$ is also a set of uniqueness for an Ising model induced by any graph $G$ having $k$~vertices. That includes, e.g., the set $W_{\{1,k\}}$ consisting of only $k+1$ elements, cf. Lemma~\ref{rem:three}. However, the set $W_{\{\floor{{k}/{2}},\floor{{k}/{2}}+1\}}$ is not of uniqueness for $\ising$\\ for any $k\geq 2$.

\section{Extremal properties}
\label{sec:extr}

To improve the probability bounds to achieve the existence of the MLE for iid samples, we discuss the sets of uniqueness for $(\set{B}^k_q)_+$ of the smallest cardinality. 
For the latter discussion, we introduce the following notation:
\begin{enumerate}
    \item $u(k,q)$:= size of the smallest set of uniqueness for $\bkqp{k}{q}$,
    \item $g(k,q)$:= size of the smallest subset of $\set{X}=\{-1,+1\}^k$ which intersects~every $(k-q)$-subcube.
\end{enumerate}

Observe that every set of uniqueness for $\bkqp{k}{q}$ intersects every {$(k-q)$-subcube}, because for any $(k-q)$-subcube $S$, the function $\varphi(x)=\one_{\{x\in S\}} \in \bkqp{k}{q}$. Thus, we obtain the following result.
\begin{corollary}\label{cor:ineq}
For every $k$ and $q$ we have 
$$
g(k,q) \leq 
u(k,q).
$$
\end{corollary}
The result above gives a lower bound for $u(k,q)$, which was studied as one of~the~T\'uran-type problems, cf.~\cite{MR1204793}. 

For larger values of $q$, we use the fact proved by Graham et al.~\cite{MR1204793} that $g(k,q)$ can~be~bounded from below with the values of $g$ for lower $q$. 

\begin{lemma}(\cite[\rm{Theorem 2. (ii)}]{MR1204793})\label{lem:lowerbound_gkq}
For $k\geq 1$ and $q\leq k-1$ we have 
$$
g(k,q) \geq 2g(k-1,q-1).
$$
\end{lemma}

The exact value of $g(k,2)$ was established by Graham et al.~\cite{MR1204793}.

\begin{lemma}(\cite[\rm{Theorem 9}]{MR1204793})
Let $k\geq 4$. Then 
$$
g(k,2) = \min \{r: \binom{r-1}{\floor{{k}/{2}}-1} \geq k \}.
$$
\end{lemma}

\begin{corollary}(\cite[\rm{Corollary 9.1}]{MR1204793})~\label{lem:uk2}
$$g(k,2) \geq \log_2 k + \frac12 \log_2 \log_2 k + O(1),$$ 
where $O(1)$ is a non-negative term.
\end{corollary}
The upper bounds for $u(k,q)$ can be established with the~size of a known set of~uniqueness either for $\bkqp{k}{q}$ or for a larger class $\bkq{k}{q}$. 
\subsection{Partial results in the general case}
Below we state some partial results for~the~more general case $(\set{B}^k_q)_+$. To simplify the further notation, we denote the following equations and inequalities referring to the level sets $W_j$, $j=0,1,\ldots,k$:
\begin{align}
0 &= \sum\limits_{x\in W_j} \varphi(x) = \binom{k-q}{j} T_0 + \ldots + (k-q) T_{j-1} + T_j\tag{equation\ $E_j$}\\
0 &\leq \sum\limits_{x\in W_j} \varphi(x) = \binom{k-q}{j} T_0 + \ldots + (k-q) T_{j-1} + T_j\tag{inequality $I_j$}
\end{align}
as $E_j$ and $I_j$, respectively.
\begin{theorem}\label{thm:uk3}
Let $k\geq 4$. Then $u(k,3)\leq{2k}$.
\end{theorem}
\begin{proof}
Specifically, we prove that the set $W_{\{1,k-1\}}$, which consists of $2k$ elements, is of uniqueness for $\bkqp{k}{3}$. Indeed, let $\varphi\in\bkqp{k}{3}$ vanish on $W_{\{1,k-1\}}$. It implies the equations $E_1$ and $E_{k-1}$ as well as the inequalities $I_j$ for every\\ $0\leq j \leq k$. In particular, we obtain the system $\{I_0, E_1, I_2, E_{k-1}, I_k\}$, i.e.,
\begin{equation*}
\begin{cases}
T_0 &\geq 0,\\
(k-3)T_0 + T_1&= 0,\\
\binom{k-3}{2}T_0 + (k-3)T_1 + T_2  &\geq 0,\\
T_2 +(k-3)T_3 &= 0,\\
T_3 &\geq 0.
\end{cases}
\end{equation*}
This yields the following chain of inequalities:
\begin{align*}
0 &\overset{I_2}{\leq} \binom{k-3}{2}T_0 + (k-3)T_1 + T_2\\ &= \frac{(k-3)(k-4)}{2}T_0 + (k-3)T_1 + T_2 \\ &\overset{E_1}{=}\left({(k-3)(k-4)}/{2}  - (k-3)^2\right) T_0 + T_2\\ 
&=  (k-3)\left({k}/2 - 2 - k + 3\right) T_0 + T_2 \\
&=  (k-3)\left(1-{k}/2\right) T_0 + T_2 \\
&\overset{E_{k-1}}{=}
(k-3) \left(1-{k}/{2}\right) T_0 - (k-3)T_3\\
&\overset{I_k}{\leq} (k-3) \left(1-{k}/{2}\right) T_0\\  &\overset{I_0}{\leq} 0.
\end{align*} 
Hence, equality holds throughout. In particular, the equations $E_0: T_0=0$ and $E_k: T_3=0$ hold. Now, the equalities $E_1$ and $E_{k-1}$ imply that $T_1=0$ and $T_2=0$, respectively. Therefore, $W_{\{1,k-1\}}$ is~a~set~of~uniqueness for $\bkqp{k}{3}$.
\end{proof}
\begin{theorem}~\label{thm:blofeld}
For any $q < k$ there exists a set of uniqueness for $\bkqp{k}{q}$ consisting of $\ 2\left(1+k+\binom{k}{2} + \ldots + \binom{k}{\lfloor{q}/{2}\rfloor}\right) = O(k^{\floor{{q}/{2}}})$ elements.
\end{theorem}
Specifically, we prove that the set $U=W_{\{0,1,2,\ldots,\floor{{q}/{2}},k-\floor{{q}/{2}},\ldots,k-1,k \}}$, consisting of $2\left(1+k+\binom{k}{2} + \ldots + \binom{k}{\lfloor{q}/{2}\rfloor}\right)$ elements, is a set of uniqueness for~the~linear space $\bkq{k}{q}$. Therefore, it is also of uniqueness for $\bkqp{k}{q}$.
\begin{proof}
Let $U=W_{\{0,1,2,\ldots,\floor{{q}/{2}},k-\floor{{q}/{2}},\ldots,k-1,k \}}$. Assume that $\varphi(U)=\{0\}$. It implies the system of equations 
$$
\{E_0,E_1,E_2,\ldots,E_{\floor{\frac{q}{2}}},E_{k-\floor{\frac{q}{2}}},\ldots,E_{k-1},E_k\},
$$
which is the following:
\begin{equation*}
\begin{cases}
T_0=0, \\
(k-q) T_0 + T_1 = 0, \\
\binom{k-q}{2}T_0 + (k-q)T_1 + T_2 = 0, \\
\cdots \\
\binom{k-q}{\floor{{q}/{2}}} T_0\ \ +\ \ \binom{k-q}{\floor{{q}/{2}} - 1}\ T_1\  + \ldots + (k-q)T_{\floor{{q}/{2}} - 1} + T_{\floor{{q}/{2}}} = 0, \\
----------------------------\\
T_{q-\floor{{q}/{2}}} + (k-q) T_{q-\floor{{q}/{2}} + 1} + \ldots + \binom{k-q}{\floor{{q}/{2}} - 1} T_{q-1} + \binom{k-q}{\floor{{q}/{2}}} T_q = 0,\\ 
\cdots \\
T_{q-1} + (k-q) T_q = 0,\\
T_q = 0.
\end{cases}
\end{equation*}
It is equivalent to the matrix equation $\mathbf{A} \mathbf{T} = \mathbf{0}$, where $\mathbf{T} = [T_0,T_1,\ldots,T_q]'$ and $\mathbf{A}$ is invertible.
Therefore the only solution of the system of equations is\\ $T_0=T_1=\ldots=T_q=0$ and~thus $\varphi(\set{X})=\{0\}$.\\
Summarizing, the set $U$, having $2\left(1+k+\binom{k}{2} + \ldots + \binom{k}{\lfloor{q}/{2}\rfloor}\right) = O(k^{\floor{{q}/{2}}})$ elements, is of uniqueness for $\bkq{k}{q}$ and thus also for $\bkqp{k}{q}$.
\end{proof}

We conclude the discussion on the smallest size of the set of uniqueness for $\bkqp{k}{q}$ with the following theorem.

\begin{theorem}
Let $k\geq 4$ and $q\leq k$. Then the following statements hold:
\begin{enumerate}
    \item $u(k,0)=1$.
    \item $u(k,1)=2$.
    \item $\log_2 k + \frac12 \log_2 \log_2 k + O(1) \leq u(k,2) \leq k+1$.
    \item $2 \log_2 (k-1) + \log_2 \log_2 (k-1) + O(1) \leq u(k,3) \leq 2k$.
    \item $u(k,k-1) = 2^{k-1}$.
    \item $u(k,k) = 2^k$.
\end{enumerate}
Moreover, for any $3\leq q\leq k-1$ we obtain
$$
2^{q-2}\log_2(k-q+2) + 2^{q-3}\log_2\log_2(k-q+2)+O(1) \leq u(k,q) \leq O\left(k^{\lfloor{q}/{2}\rfloor} \right).
$$
\end{theorem}
\begin{proof}
(1): 
When $q=0$, all functions in $\set{B}=(\set{B}^k_0)_+$ are the constant positive functions. Therefore, every non-empty subset of $\set{X}$ is of uniqueness, thus for every $k$ the smallest size equals $1$.\\
(2): 
When $q=1$, by Remark~\ref{rem:012}, the smallest set of uniqueness is $\{x,-x\}$, $x\in\set{X}$, and its size equals $2$.\\
(3): 
By Lemma~\ref{rem:three}, the set $W_{\{1,k\}}$ is a minimal set of uniqueness for $(\set{B}^k_2)_+$, therefore the smallest cardinality of the set of uniqueness does not exceed $k+1$. The lower bound is the result of Corollary~\ref{lem:uk2}.\\ 
(4):
The upper bound for $u(k,3)$ follows directly from Theorem~\ref{thm:uk3}. The lower bound is derived from Corollary~\ref{cor:ineq}, Lemma~\ref{lem:lowerbound_gkq} and Corollary~\ref{lem:uk2}. Indeed,
\begin{align*}
u(k,3)&\overset{\rm{Cor.}~\ref{cor:ineq}}{\geq}g(k,3)\\
&\overset{\rm{Lemma}~\ref{lem:lowerbound_gkq}}{\geq}2g(k-1,2)\\
&\overset{\rm{Cor.}~\ref{lem:uk2}}{\geq} 2\log_2(k-1)+\log_2\log_2(k-1)+O(1).
\end{align*}
(5): 
For every $k\geq 1$, the only minimal sets of uniqueness are $W_{\{0,2,4,\ldots\}}$ and $W_{\{1,3,5,\ldots\}}$, cf. Remark~\ref{rem:012}.
Both have size $2^{k-1}$.\\
(6):
By Remark~\ref{rem:012}, the only set of uniqueness consists of $2^k$ elements.\\
(General case):
The upper bound for $u(k,q)$ is obtained by Theorem~\ref{thm:blofeld}. Similarly, as~in~part~(4) of this proof, the lower bound comes from Corollary~\ref{cor:ineq}, Lemma~\ref{lem:lowerbound_gkq} and Corollary~\ref{lem:uk2}. Indeed,
\begin{align*}
u(k,q)&\overset{\rm{Cor.}~\ref{cor:ineq}}{\geq}g(k,q)\\
&\overset{\rm{Lemma}~\ref{lem:lowerbound_gkq}}{\geq \ldots \geq}
2^{q-2}g(k-q+2,2)\\
&\overset{\rm{Cor.}~\ref{lem:uk2}}{\geq}2^{q-2}\log_2(k-q+2)+2^{q-3}\log_2\log_2(k-q+2)+O(1).
\end{align*}
\end{proof}

Generally speaking, the upper bounds are constructed as the sizes of the specific sets of uniqueness for $\bkqp{k}{q}$. The lower bounds are also the lower bounds for the more relaxed conditions of intersecting all $(k-q)-$subcubes. Therefore, it should be possible to improve the existing bounds. That, if possible, may result in obtaining smaller sets of uniqueness for $\bkqp{k}{q}$, both when the sets are fixed and when they are randomly sampled from $\set{X}$.

We finish this section with the following conjecture on the smallest size of the set~of~uniqueness for the Ising model.
\begin{conjecture}
For every $k\geq 3$ we have $u(k,2)=k+1$, i.e., there are no sets of~uniqueness having at most $k$ elements.
\end{conjecture}

\section{Proofs}\label{sec:proofs}
\subsection{Main characterization of the sets of uniqueness for $\ising$}\label{sec:proof:main}
\begin{proof}{Proof of Theorem {\ref{thm:main}}}
Step 1:\\
We commence the proof by mapping each single level set $W_j$: $j=0,1,\ldots,k$ to~a~point 
$$
V^{(j)}:=(V_0^{(j)}, V_1^{(j)}, V_2^{(j)}) = \left(\binom{k-2}{j},\binom{k-2}{j-1},\binom{k-2}{j-2}\right)\in\R^3.
$$
The value of $V_i^{(j)}$ is equal to $|W_j \cap S|$ for $S\in S_i \subset \set{S}^k_2$, cf. Lemma~\ref{lem:card:weights}.\\
Next, we normalize each $V^{(j)}$ so that its coordinates sum to $1$. Therefore, we~perform the following transformation:
\begin{align*}
h&:\ \R^3 \to \mathbb{D}^2 \subset \R^3,\\
h&:\ V^{(j)} \mapsto Q^{(j)} = (Q_1^{(j)},Q_2^{(j)},Q_3^{(j)}) := {V^{(j)}}/({V_0^{(j)}+V_1^{(j)}+V_2^{(j)}}),
\end{align*}
where the hyperplane $\mathbb{D}^2 := \{(v_1,v_2,v_3)\in\R^3:\ v_1+v_2+v_3=1\}$. Note that within $\R^3_+$ this hyperplane covers one of the facets of~the~unit ball in~the~$\ell_1$-norm, i.e., $\mathbb{D}^2\ \cap\ \R^3_+\ = B_1(0,1)\ \cap\ \R^3_+$.
Now we project the set $\{Q^{(0)},Q^{(1)},\ldots,Q^{(k)}\}$ onto $\R^2$, where each $Q^{(j)} = (Q_1^{(j)}, Q_2^{(j)}, Q_3^{(j)})$ is mapped to 
$$P^{(j)}=(P_1^{(j)},P_2^{(j)}):=(Q_1^{(j)},Q_2^{(j)}).$$
Since projections preserve convexity, cf.~\cite[Sec. 2.3.2.]{boyd}, the polygonal chain\\ $[Q^{(0)},Q^{(1)},\ldots,Q^{(k)},Q^{(0)}]$ bounds a~convex polytope in $\R^3$ if and only if\\ $[P^{(0)},P^{(1)},\ldots,P^{(k)},P^{(0)}]$ bounds a~convex polytope in $\R^2$.\\
Step 2:\\
Now we want to prove that $[P^{(0)},P^{(1)},\ldots,P^{(k)},P^{(0)}]\subset\R^2$ indeed bounds a~convex polytope in $\R^2$. We start by observing the monotonicity of the slope\\ of~the~points $P^{(j)}$. For that reason, we introduce the following lemmas, which we~prove in Section~\ref{sec:other}: 
\begin{lemma}\label{lem:symmetry}
Let $0\leq j\leq k$. Then $P_2^{(j)} = P_{2}^{(k-j)}$.
\end{lemma}
\begin{figure}
    \centering
    \begin{tikzpicture}[line cap=round,line join=round,x=1cm,y=1cm]
\begin{axis}[
x=1cm,y=1cm,
axis lines=middle,
xmin=-1.7223021279184858,
xmax=7.575644861458564,
ymin=-1.2045394667942726,
ymax=5.467108145173295,
xtick={0},
ytick={0},
xticklabel=\empty,yticklabel=\empty,]
\clip(-1.7223021279184858,-1.2045394667942726) rectangle (7.575644861458564,5.467108145173295);
\draw [line width=2.8pt,color=qqqqff] (0,0)-- (0,4);
\draw [line width=2.8pt,color=qqqqff] (0,0)-- (5,0);
\draw [line width=2.8pt,color=qqqqff] (5,0)-- (4.8,4);
\draw [line width=1pt] (5,-1.2045394667942726) -- (5,5.467108145173295);
\draw [line width=0.5pt,color=black] (-1.7223021279184858,4.390941218137013)-- (7.575644861458564,4.390941218137013);
\draw [line width=0.5pt,color=black] (-1.7223021279184858,5)-- (7.575644861458564,5);
\draw [line width=0.5pt,color=black] (-1.7223021279184858,3.447613328493216)-- (7.575644861458564,3.447613328493216);
\draw [line width=0.5pt,color=black] (-1.7223021279184858,4.6509085559997935)-- (7.575644861458564,4.6509085559997935);

\draw [line width=2.8pt,color=qqqqff] (0,4)-- (0.6858315573014054,4.390941218137013)-- (1,5)-- (1.3517288539867929,3.447613328493216)-- 
(2.2584826196860437,4.6509085559997935)-- 
(3.717789461358276,3.447613328493216)-- (4.199502399386003,5)-- (4.52537,4.390941218137013)-- (4.8,4);
\draw [color=qqqqff](0.1415232366946812,0.6063363135969243) node[anchor=north west] {$P^{(k)} = (0,0)$};
\draw [color=qqqqff](5.118784153559389,0.5957463967525314) node[anchor=north west] {$P^{(0)} = (1,0)$};
\draw [color=qqqqff](4.96114382093696,4.143368539624174) node[anchor=north west] {$P^{(1)}$};
\draw [color=qqqqff](-0.03562406825075127,4.143368539624174) node[anchor=north west] {$P^{(k-1)}$};
\begin{scriptsize}
\draw [fill=ududff] (5,0) circle (2.5pt);
\draw [fill=ududff] (0,0) circle (2.5pt);
\draw [fill=ududff] (0,4) circle (2.5pt);
\draw [fill=ududff] (4.8,4) circle (2.5pt);
\draw [fill=ududff] (0.6858315573014054,4.390941218137013) circle (2.5pt);
\draw [fill=ududff] (1,5) circle (2.5pt);
\draw [fill=ududff] (1.3517288539867929,3.447613328493216) circle (2.5pt);
\draw [fill=ududff] (2.2584826196860437,4.6509085559997935) circle (2.5pt);
\draw [fill=ududff] (3.717789461358276,3.447613328493216) circle (2.5pt);
\draw [fill=ududff] (4.199502399386003,5) circle (2.5pt);
\draw [fill=ududff] (4.52537,4.390941218137013) circle (2.5pt);
\end{scriptsize}
\end{axis}
\end{tikzpicture}
\caption{Illustration of Lemma~\ref{lem:symmetry}, i.e., $P_2^{(j)} = P_{2}^{(k-j)}$.}
\label{fig1}
\end{figure}
\begin{lemma}\label{lem:slope}
Let $1\leq j\leq k-2$. Then ${P_2^{(j)}}/{P_1^{(j)}} = {j}/(k-j-1).$
\end{lemma}
\begin{corollary}\label{cor:incr}
Let $j\in\{0,1,\ldots,k-2\}$. Then the slope ${P_2^{(j)}}/{P_1^{(j)}}$ is an increasing function of $j$.
\end{corollary}
\begin{proof}
Indeed, 
$$
{P_2^{(0)}}/{P_1^{(0)}} = 0/1 = 0 < {P_2^{(1)}}/{P_1^{(1)}}.
$$ 
For larger values of $j$ we apply Lemma~\ref{lem:slope}.
\end{proof}
\begin{figure}
    \centering
    \begin{tikzpicture}[line cap=round,line join=round,x=1cm,y=1cm]
\begin{axis}[
x=1cm,y=1cm,
axis lines=middle,
xmin=-1.7223021279184858,
xmax=7.575644861458564,
ymin=-1.2045394667942726,
ymax=5.467108145173295,
xtick={0},
ytick={0},
xticklabel=\empty,yticklabel=\empty,]
\clip(-1.7223021279184858,-1.2045394667942726) rectangle (7.575644861458564,5.467108145173295);
\draw [line width=1pt] (5,-1.2045394667942726) -- (5,5.467108145173295);
\draw [color=qqqqff](-0.88010079488109762,-0.05427968125305666) node[anchor=north west] {$P^{(k)} = (0,0)$};
\draw [color=qqqqff](5.118784153559389,0.5957463967525314) node[anchor=north west] {$P^{(0)} = (1,0)$};
\draw [color=qqqqff](4.939963987248175,3.12749953583839515) node[anchor=north west] {$P^{(1)}$};
\draw [color=qqqqff](-1.39938717739688985,3.12749953583839515) node[anchor=north west] {$P^{(k-1)}$};
\draw [line width=0.7pt,domain=0:7.575644861458564] plot(\x,{(-0--3*\x)/4.8});
\draw [line width=0.7pt,domain=0:7.575644861458564] plot(\x,{(-0--4.429296294422785*\x)/2.994651398692152});
\draw [line width=0.7pt,domain=0:7.575644861458564] plot(\x,{(-0--4.429296294422785*\x)/2.005348601307848});
\draw [line width=0.7pt,domain=0:7.575644861458564] plot(\x,{(-0--4*\x)/4});
\draw [line width=0.7pt,domain=0:7.575644861458564] plot(\x,{(-0--4*\x)/1});
\draw [color=qqqqff](3.9316338727117774,4.030361782692849) node[anchor=north west] {$P^{(2)}$};
\begin{scriptsize}
\draw [fill=xdxdff] (5,0) circle (2.5pt);
\draw [fill=ududff] (0,0) circle (2.5pt);
\draw [fill=ududff] (0,3) circle (2.5pt);
\draw [fill=ududff] (4.8,3) circle (2.5pt);
\draw [fill=ududff] (1,4) circle (2.5pt);
\draw [fill=ududff] (2.005348601307848,4.429296294422785) circle (2.5pt);
\draw [fill=ududff] (4,4) circle (2.5pt);
\draw [fill=ududff] (2.994651398692152,4.429296294422785) circle (2.5pt);
\end{scriptsize}
\end{axis}
\end{tikzpicture}
\caption{Illustration of the Corollary~\ref{cor:incr}.}
\label{fig2}
\end{figure}
\begin{lemma}\label{lem:10}
For $j<\frac{k}2$ the value of $P_2^{(j)}$ increases with $j$ and for~$j\geq\frac{k}2$ it~decreases with $j$.
\end{lemma}
By the definition of $P^{(j)}$ and Lemma~\ref{lem:card:weights}, we obtain $P^{(0)}=(1,0)$ and $P^{(k)}=(0,0)$. Moreover, according to~Lemma~\ref{lem:symmetry}, the second coordinates of $P^{(j)}$ are distributed symmetrically, i.e., $P_2^{(j)} = P_{2}^{(k-j)}$.\\
By Corollary~\ref{cor:incr}, the slope ${P_2^{(j)}}/{P_1^{(j)}}$ is an increasing function of~$j$.\\ 
Moreover, 
\begin{align*}
(P_{1}^{(k-1)},P_{2}^{(k-1)}) = (0,{1}/({k-1}))&& \textrm{and} &&(P_1^{(k)},P_2^{(k)}) = (0,0).
\end{align*}
By Lemma~\ref{lem:10} the value of $P_2^{(j)}$ increases with $j$ for $j<{k}/2$ and decreases with~$j$ for ${j\geq{k}/2}$, see Figure~\ref{fig2}.\\
Therefore, the polygonal chain $[P^{(0)},P^{(1)},\ldots,P^{(k)},P^{(0)}]$ is not self-intersecting.\\ 
Step 3:\\
To~prove that $[P^{(0)},P^{(1)},\ldots,P^{(k)},P^{(0)}]$ bounds a convex polytope, we want \\ to show that no point $P^{(j)}$ lies~in~the~interior of~${\rm{Conv}}(\{P^{(0)},P^{(1)},\ldots,P^{(k)}\})$.\\ Given $T\in\R^3$, defined as in~\eqref{eq:Tdef}, we introduce the following map $\psi$:
\begin{align}~\label{eq:psi}
\psi&:\ \mathbb{D}^2\to \R\\
\psi&:\ Q^{(j)} \mapsto Q_1^{(j)} T_0 + Q_2^{(j)} T_1 + Q_3^{(j)} T_2 = \left\langle Q^{(j)}, T \right\rangle.
\nonumber
\end{align}
Note that $\psi$, as an inner product, is a linear map, that is
\begin{align*}
\psi\left(\lambda Q^{(j_1)} + (1-\lambda) Q^{(j_2)}\right) = \lambda \psi(Q^{(j_1)}) + (1-\lambda)\psi(Q^{(j_2)}).
\end{align*}
Now we are ready to connect the map $\psi$ with the sets of uniqueness.\\ 
Observe that
\begin{align*}
\sum\limits_{x\in W_j}\varphi(x) &= \left\langle V^{(j)}, T\right\rangle\\
&= \left(V_0^{(j)}+V_1^{(j)}+V_2^{(j)}\right)\left\langle Q^{(j)}, T\right\rangle\\
&= \left(V_0^{(j)}+V_1^{(j)}+V_2^{(j)}\right)\psi\left(Q^{(j)}\right).
\end{align*}
Since $\varphi$ is non-negative,
\begin{align*}
\sum\limits_{x\in W_j}\varphi(x) = 0 && {\rm{if\ and\ only\ if}} && \varphi(W_j) = \{0\}.
\end{align*}
Therefore, we obtain the following proposition.
\begin{proposition}~\label{lem:psizero}
Let $k\geq 3$, $0 \leq j \leq k$. Let $\varphi\in\ising$. Then 
\begin{align*}
\varphi(W_j)=\{0\} && {\rm{if\ and\ only\ if}} && \psi\left(Q^{(j)}\right)=0.
\end{align*}
\end{proposition}
Suppose that $[P^{(0)},P^{(1)},\ldots,P^{(k)},P^{(0)}]$ does not bound a convex polytope, i.e., there exists $0\leq j\leq k$ such that $P^{(j)}\in{\rm{int}}\left({\rm{Conv}}(\{P^{(0)},P^{(1)},\ldots,P^{(k)}\})\right)$. Then also\\ $Q^{(j)} \in {\rm{int}}\left({\rm{Conv}}(\{Q^{(0)},Q^{(1)},\ldots,Q^{(k)}\})\right)$, which implies the existence of 
\begin{align*}
\lambda_0,\lambda_1,\ldots,\lambda_{j-1},\lambda_{j+1},\ldots,\lambda_k>0& &\rm{such\ that}& &Q^{(j)} = \sum\limits_{\substack{i=0\\ i\neq j}}^{k} \lambda_i Q^{(i)}.
\end{align*} 
Therefore, 
$$
\psi(Q^{(j)}) = \sum\limits_{\substack{i=0\\ i\neq j}}^{k} \lambda_i \psi(Q^{(i)})
$$ 
with $\psi$ defined in~\eqref{eq:psi}.\\
Now assume that $\varphi(W_j)=\{0\}$. Then 
\begin{align*}
\sum_{x\in W_j}\varphi(x) = 0&& \rm{and}&& \psi(Q^{(j)}) = 0.
\end{align*}
Thus also $\psi(Q^{(i)}) = 0$ for every $i\neq j$. But then $\varphi(\set{X})=\{0\}$, so $W_j$ is a set of uniqueness for $\ising$.\\
That way we obtain the following implication:
\begin{align}\label{eq:1el}
&\left(\forall\ 0\le j\le k, \ W_j {\rm{\ is\ not\ of\ uniqueness\ for\ }}\ising\right)\Then \\ &\Then \left(\left[P^{(0)},P^{(1)},\ldots,P^{(k)},P^{(0)}\right]{\rm{\ bounds\ a\ convex\ polytope}}\right).
\nonumber
\end{align}
Now we proceed with proving the implication:
\begin{align}\label{eq:pol}
&\left(\left[P^{(0)},P^{(1)},\ldots,P^{(k)},P^{(0)}\right]{\rm{\ bounds\ a\ convex\ polytope}}\right) \Then \\ &\Then \left(W_{\{j,j+1\}} {\rm{\ is\ not\ of\ uniqueness\ for\ }}\ising\ \forall\ j=0,1,\ldots,k-1\right).
\nonumber
\end{align}

For $j=0$, it suffices to consider a function $\varphi=\one_{\{x\in(+1,+1,*,*,\ldots,*)\}}\in\ising$, which vanishes on $U=W_{\{0,1\}}$. An analogous argument holds for $j=k-1$ with a function $\varphi=\one_{\{x\in(-1,-1,*,*,\ldots,*)\}}\in\ising$ vanishing on~$U=W_{\{k-1,k\}}$.\\ 
Now, let us consider $U=W_{\{j,j+1\}}$ for $1\leq j \leq k-2$. We seek for a nontrivial function $\varphi\in\ising$, satisfying $\varphi(U)=\{0\}$.
Under the assumption that $P$ (and $Q$) bounds a convex polytope, we~show the existence of a vector $T=(T_0,T_1,T_2)$ such that
$$\psi(Q^{(j)})=\psi(Q^{(j+1)})=0$$
and
\begin{align*}
\psi(Q^{(i)})>0, &&\mathrm{for\ all}\ {0\leq i \leq k, i\neq j, j+1}. 
\end{align*} 
Then, thanks to Proposition~\ref{lem:psizero}, by setting
\begin{align*}
\varphi(x)=\frac{1}{|W_i|}\psi(Q^{(i)})\geq 0, &&\mathrm{for\ all}\ x\in W_i, {0\leq i \leq k,} 
\end{align*} 
we conclude that
$$
\sum\limits_{x\in W_j}\varphi(x) = \sum\limits_{x\in W_{j+1}}\varphi(x) = 0,
$$
and
\begin{align*}
\sum\limits_{x\in W_i}\varphi(x)>0, &&\mathrm{for\ all}\ 0\leq i \leq k, i\neq j, j+1.
\end{align*} 
Now our problem is reduced to finding $(T_0,T_1,T_2)$. Observe that $[Q^{(j)},Q^{(j+1)}]$ is an edge of ${\rm{Conv}}(Q)\subset\mathbb{D}^2$. Since ${\rm{Conv}}(Q)$ is a convex polytope, for each $z\in{\rm{Conv}}(Q)$ there are $\lambda_1,\lambda_2\in\R$ and $\lambda_2\geq0$ such that
\begin{equation}\label{eq:sum}
z = Q^{(j)} + \lambda_1 \left( Q^{(j+1)} - Q^{(j)}\right) + \lambda_2 \left( Q^{(0)} - Q^{(j)}\right).
\end{equation}
Let $\bf{s}=\pm Q^{(j+1)}\times Q^{(j)}$ be the cross product of $Q^{(j+1)}$ and $Q^{(j)}$, where the sign is chosen in such a way that
$$
\langle {\bf{s}}, Q^{(0)} \rangle < 0.
$$
Observe that
$$
\langle {\bf{s}}, Q^{(j)} \rangle = \langle {\bf{s}}, Q^{(j+1)} \rangle = 0.
$$
Now, for each $z\in{\rm{Conv}}(Q)$, by~\eqref{eq:sum}, we get
\begin{align*}
\langle {\bf{s}},{z}\rangle &= \langle {\bf{s}},Q^{(j)}\rangle + \lambda_1 \langle {\bf{s}},Q^{(j+1)} - Q^{(j)}\rangle + \lambda_2  \langle {\bf{s}},Q^{(0)}-Q^{(j)}\rangle\\
&= \lambda_2 \langle {\bf{s}},Q^{(0)}\rangle \leq 0.
\end{align*}
Moreover, $\langle {\bf{s}},{z}\rangle = 0$ for $z\in{\rm{Conv(Q)}}$, if and only if $\lambda_2 = 0$, that is
$$
\bf{z} \in [Q^{(j)}, Q^{(j+1)}].
$$
Hence, it is enough to take $(T_0,T_1,T_2) = -{\bf{s}}$. Indeed, then we have
$$
\psi\left(Q^{(j)}\right) = -\langle\bf{s},Q^{(j)}\rangle = 0,
$$
and
$$
\psi\left(Q^{(j+1)}\right) = -\langle\bf{s},Q^{(j+1)}\rangle = 0,
$$
and
$$
\psi\left(Q^{(i)}\right) = -\langle\bf{s},Q^{(i)}\rangle > 0
$$
for all $0\leq i \leq k$, $i\neq j, j+1$. This completes the proof of~\eqref{eq:pol}.
\\
Take $(T_0,T_1,T_2)={\bf{T}}={\bf{-s}}$. By Proposition~\ref{lem:psizero} and the properties of $\psi$, this vector results in a function $\varphi\in\ising$ vanishing on $W_j$ and $W_{j+1}$. Indeed, observe that
\begin{equation*}
\begin{cases}
\psi(Q^{(j)}) = - \langle {\bf{s}},Q^{(j)}\rangle &= 0,\\
\psi(Q^{(j+1)}) = - \langle {\bf{s}},Q^{(j+1)}\rangle &= 0.
\end{cases}
\end{equation*}
Additionally, for every $i\neq j,j+1$ we have $\varphi(W_i)> 0$, since ${-\langle {\bf{s}},Q^{(i)}\rangle > 0}$.
This completes the proof of the implication~\eqref{eq:pol}.\\
Observe that for every $j=0,1,\ldots,k-1$ we have the implication
\begin{align}\label{eq:obvious}
\left(W_{\{j,j+1\}} {\rm{\ is\ not\ of\ uniqueness\ for\ }}\ising\ \forall\ j=0,1,\ldots,k-1\right) 
\Then\\
\Then\left(W_{j}, W_{j+1} {\rm{\ are\ not\ of\ uniqueness\ for\ }}\ising\ \forall\ j=0,1,\ldots,k-1\right).
\nonumber
\end{align}

From equations~\eqref{eq:1el},~\eqref{eq:pol} and~\eqref{eq:obvious}, for every $j=0,1,\ldots,k-1$ we obtain that $W_j$ is of uniqueness for $\ising$ if and only if $W_{j+1}$ is of uniqueness for $\ising$, too. But $W_0$ is~not of~uniqueness, neither is any other $W_j$.\\ 
Step 4:\\
We finish the proof by showing that if $2\leq j-i\leq k-1$, then $W_{\{i,j\}}$ is a set of uniqueness for $\ising$.
First we observe that $U=W_{\{0,k\}}$ is not of uniqueness for $\bkqp{k}{2}$. To realize that, one may consider the function 
$$
\varphi=\one_{\{x\in(+1,-1,*,*,\ldots,*)\}}\in\ising,
$$ 
which equals one for every point of $\set{X}$ having its first two coordinates fixed as~$+1$ and $-1$, respectively, 
and equals zero everywhere else.
Note that $\varphi$ vanishes on both $W_0=\{(-1,-1,\ldots,-1)\}$ and $W_k=\{(+1,+1,\ldots,+1)\}$.
Without loss of generality, assume that $i<j$. Since $P$ bounds a convex polytope, the segment $[P^{(i)},P^{(j)}]$ is one of its proper diagonals and~intersects every other diagonal\\ $[P^{(h_{in})},P^{(h_{ex})}]$, where $i<h_{in}<j$ and $h_{ex}\in\{0,1,\ldots,i-1,j+1,\ldots,k\}$.\\ The final statement comes from the following Lemma, which is proven in Section~\ref{sec:other}:
\begin{lemma}\label{lem:dziobak}
Let $j_1,j_2,j_3,j_4$ be the vertices of\ ${\rm{Conv}}(\{P^{(0)},P^{(1)},\ldots,P^{(k)}\})$\\ and~let~$\varphi\in\ising$. Assume that the segments $[P^{(j_1)},P^{(j_3)}]$ and $[P^{(j_2)},P^{(j_4)}]$ intersect. Then~$\varphi(W_{\{j_1,j_3\}}) = \{0\}$ if and only if $\varphi(W_{\{j_2,j_4\}}) = \{0\}$.
\end{lemma}

Indeed, by Lemma~\ref{lem:dziobak}, $\varphi\left(W_{\{h_{in},h_{ex}\}}\right)=\{0\}$. Then, by the arbitrariness of $h_{in}$ and $h_{ex}$, $\varphi\left(\set{X}\right)=\{0\}$, which implies that $W_{\{i,j\}}$ is a set of uniqueness for $\ising$ and completes the~proof.
\end{proof}

\subsection{Other proofs}
\label{sec:other}
\begin{proof}{Proof of Remark~{\ref{rem:compare}}}
Below we present the examples of functions,\\ which distinguish the Ising models on different graphs.\\
1. The set $U_1=\{(-1,-1,-1),(+1,+1,+1)\}$ is of uniqueness for $\bkqp{3}{1}$ by~\cite[Example 3.7]{MLE}. Observe that the function $\one_{\{x\in(+1,-1,*)\}}$ belongs to $\set{B}'_+$ as the characteristic function of the edge $(+1,-1,*)$, which is parallel to the edge $(-1,-1,*)$, see Example~\ref{ex:ising}. This function vanishes on $U_1$, which therefore is not of uniqueness for $\set{B}'_+$.\\ 
2. Let $U_2 = \{(+1,-1,-1), (-1,+1,-1), (-1,-1,+1), (-1,+1,+1)\}$.\\ To prove that $U_2$ is not of uniqueness for $\set{B}''_+$, we show that the function\\ $\one_{\{x\in(+1,*,+1)\}}$ belongs to $\set{B}''_+$ and vanishes on $U_2$.\\
Indeed, this function equals zero on every element of $U_2$. It is also the characteristic function of the edge $(+1,*,+1)$, which is parallel to the edge $(-1,*,-1)$. Therefore this function belongs to $\set{B}''_+$.\\
Now we show that $U_2$ is of uniqueness for $\set{B}'_+$. Indeed, let $\varphi$ be any function from $\set{B}'_+$. Then $\varphi$ is non-negative and can be represented as
\begin{align*}
\varphi &= \alpha\one_{\{x\in(-1,-1,*)\}} + \beta\one_{\{x\in(-1,+1,*)\}} + \gamma\one_{\{x\in(+1,-1,*)\}} + \delta\one_{\{x\in(+1,+1,*)\}}\\ &+ A\one_{\{x\in(*,*,-1)\}} + B\one_{\{x\in(*,*,+1)\}}.
\end{align*}
Suppose that $\varphi(U_2)=\{0\}$. Then we have $\alpha + \gamma + A + B = 0$, because
\begin{equation*}
    \begin{cases}
        0 = \varphi((-1,-1,+1)) = \alpha + B,\\
        0 = \varphi((+1,-1,-1)) = \gamma + A.
    \end{cases}
\end{equation*}
Since $\varphi$ is non-negative, we also have
\begin{equation*}
    \begin{cases}
        0 \leq \varphi((-1,-1,-1)) = \alpha + A,\\
        0 \leq \varphi((+1,-1,+1)) = \gamma + B,
    \end{cases}
\end{equation*}
which implies that
\begin{equation*}
\varphi((-1,-1,-1)) = \varphi((+1,-1,+1))=0.
\end{equation*}
Analogously, one can deduce that $\varphi((-1,+1,+1)) = \varphi((+1,+1,+1))=0$,\\ which implies that $\varphi(\set{X})=\{0\}$, so $U_2$ is of uniqueness for $\set{B}'_+$.\\
3. The set $U_3=\set{X}\setminus\{(-1,-1,-1),(+1,+1,+1)\}$ is not of uniqueness for $\bkqp{3}{2}$ by~\cite[Remark 5.4]{MLE}. Now we show that $U_3$ is of uniqueness for $\set{B}''_+$.\\ Indeed, let $\varphi\in\set{B}''_+$. Then $\varphi\geq 0$ and can be represented as
\begin{align*}
\varphi &= \alpha\one_{\{x\in(-1,-1,*)\}} + \beta\one_{\{x\in(-1,+1,*)\}} + \gamma\one_{\{x\in(+1,-1,*)\}} + \delta\one_{\{x\in(+1,+1,*)\}}\\ &+ \varepsilon\one_{\{x\in(-1,*,-1)\}} + \zeta\one_{\{x\in(-1,*,+1)\}} + \eta\one_{\{x\in(+1,*,-1)\}} + \theta\one_{\{x\in(+1,*,+1)\}}.
\end{align*}
Suppose that $\varphi(U_3)=\{0\}$. Then
\begin{align*}
&\varphi((-1,-1,-1)) = \alpha + \varepsilon = (\alpha + \zeta) - (\beta + \zeta) + (\beta + \varepsilon)\\ = &\varphi((-1,-1,+1)) - \varphi((-1,+1,+1)) + \varphi((-1,+1,-1)) = 0
\end{align*}
and
\begin{align*}
&\varphi((+1,+1,+1)) = \delta + \theta = (\delta + \eta) - (\gamma + \eta) + (\gamma + \theta)\\ 
= &\varphi((+1,+1,-1)) - \varphi((+1,-1,-1)) + \varphi((+1,-1,+1)) = 0.
\end{align*}
Therefore, $\varphi(\set{X})=\{0\}$, so $U_3$ is of uniqueness for $\set{B}''_+$.
\end{proof}

\begin{proof}{Proof of Lemma~{\ref{lem:symmetry}}}
By the definition of $P^{(j)}$ and Lemma~\ref{lem:card:weights}, we have
\begin{align*}
P^{(j)} = \frac{\left( \binom{k-2}{j}, \binom{k-2}{j-1}, \binom{k-2}{j-2} \right)}{\binom{k-2}{j} + \binom{k-2}{j-1} + \binom{k-2}{j-2}}
&&
and
&&
P^{(k-j)} = \frac{\left( \binom{k-2}{k-j}, \binom{k-2}{k-j-1}, \binom{k-2}{k-j-2} \right)}{\binom{k-2}{k-j} + \binom{k-2}{k-j-1} + \binom{k-2}{k-j-2}}. 
\end{align*}
Therefore, since $\binom{n}{n-m} = \binom{n}{m}$ for every $n,m\geq 0$, we obtain
\begin{align*}
P^{(j)}_2 = \frac{\binom{k-2}{j-1}}{\binom{k-2}{j} + \binom{k-2}{j-1} + \binom{k-2}{j-2}} = \frac{\binom{k-2}{k-j-1}}{\binom{k-2}{k-j} + \binom{k-2}{k-j-1} + \binom{k-2}{k-j-2}} = P^{(k-j)}_2.
\end{align*}
\end{proof}

\begin{proof}{Proof of Lemma~{\ref{lem:slope}}}
By Lemma~\ref{lem:card:weights} we~have
\begin{align*}
\frac{P_2^{(j)}}{P_1^{(j)}} &= \frac{\binom{k-2}{j-1}}{\binom{k-2}{j} + \binom{k-2}{j-1} + \binom{k-2}{j-2}}\  \cdot\  \frac{\binom{k-2}{j} + \binom{k-2}{j-1} + \binom{k-2}{j-2}}{\binom{k-2}{j}}\\ &= \frac{\binom{k-2}{j-1}}{\binom{k-2}{j}}= \frac{(k-2)!}{(j-1)! (k-j-1)!}\ \cdot\ \frac{j!(k-j-2)!}{(k-2)!}\\ 
&= \frac{j}{k-j-1}.
\end{align*}
\end{proof}

\begin{proof}{Proof of Lemma~{\ref{lem:10}}}
The proof is computational.\\ Recall that 
\begin{align*}
P_2^{(0)}=0&& {\rm{and}}&& P_2^{(1)} = \frac{\binom{k-2}{0}}{\binom{k-2}{1} + \binom{k-2}{0} + \binom{k-2}{-1}} = \frac{1}{k-1} > 0.
\end{align*}
For larger $j<\frac{k}{2}$ we will show that the ratio ${P_{2}^{(j+1)}}/{P_{2}^{(j)}}$ is greater than $1$. Indeed, let $j\geq 1$. Recall that: $$\frac{\binom{k-2}{j}}{\binom{k-2}{j-1}} = \frac{k-j-1}{j},$$ for details see the proof of Lemma~\ref{lem:card:weights}. Then
\begin{align*}
\frac{P_2^{(j+1)}}{P_{2}^{(j)}} &= \frac{\binom{k-2}{j}}{\binom{k-2}{j+1} + \binom{k-2}{j} + \binom{k-2}{j-1}}\  \cdot\  \frac{\binom{k-2}{j} + \binom{k-2}{j-1} + \binom{k-2}{j-2}}{\binom{k-2}{j-1}}\\ &= \frac{k-j-1}{j}\ \cdot\ \frac{\binom{k-2}{j} + \binom{k-2}{j-1} + \binom{k-2}{j-2}}{\binom{k-2}{j+1} + \binom{k-2}{j} + \binom{k-2}{j-1}}.
\end{align*}
Applying the same formula for different values of $k$ and $j$ yields:
\begin{align*}
\frac{P_2^{(j+1)}}{P_{2}^{(j)}}&= \frac{k-j-1}{j}\ \cdot\ \frac{\binom{k-2}{j-1}\ \cdot\ \left(\frac{k-j-1}{j} + 1 + \frac{j-1}{k-j}\right)}{\binom{k-2}{j-1}\ \cdot\ \left( \frac{(k-j-1)(k-j-2)}{j(j+1)} + \frac{k-j-1}{j} + 1 \right)}\\ &= \frac{k-j-1}{j}\ \cdot\ \frac{\frac{(k-j)(k-j-1)+(k-j)j+(j-1)j}{j(k-j)}}{\frac{(k-j-1)(k-j-2)+(k-j-1)(j+1)+j(j+1)}{j(j+1)}}\\
&=\frac{k-j-1}{k-j}\ \cdot\ \frac{j+1}{j}\\ &\cdot\ \frac{(k-j)(k-j-1)+(k-j)j+(j-1)j}{(k-j-1)(k-j-2)+(k-j-1)(j+1)+j(j+1)}\\
&=\frac{k-j-1}{k-j}\ \cdot\ \frac{j+1}{j}\ \cdot\ \frac{(k-j)(k-1)+(j-1)j}{(k-j-1)(k-1)+j(j+1)} = (*).
\end{align*}
Expanding the last fraction and straightforward calculations yields:
\begin{align*}
&\frac{(k-j)(k-1)+(j-1)j}{(k-j-1)(k-1)+j(j+1)}\\ &= \frac{(k-j)(k-1)+j(k-1)-j(k-j)}{(k-j-1)(k-1)+(j+1)(k-1)-(j+1)(k-j-1)}\\
&= \frac{k(k-1)-j(k-j)}{k(k-1)-(j+1)(k-j-1)}\\ 
&= 1 + \frac{(j+1)(k-j-1) - j(k-j)}{k(k-1) - (j+1)(k-j-1)}\\ &= 1 + \frac{(k-j)-(j+1)}{k(k-1)+(j+1)(k-j-1)}. 
\end{align*}
Therefore, for $j < {k}/{2}$ we have
\begin{align*}
(*) &= \frac{1+1/j}{1+{1}/({k-j-1})}\ \cdot\ \left(1 + \frac{(k-j)-(j-1)}{k(k-1)+(j+1)(k-j-1)}\right) \\ & > \frac{1+1/j}{1+{1}/{(k-j-1)}} > 1.
\end{align*}
Analogously, for $j \geq {k}/2$, the value of $P_2^{(j)}$ decreases with $j$ by Lemma~\ref{lem:symmetry}.
\end{proof}

\begin{proof}{Proof of Lemma{~\ref{lem:dziobak}}}
Assume that the segments $[P^{(j_1)},P^{(j_3)}]$\\ and $[P^{(j_2)},P^{(j_4)}]$ intersect. Then also $[Q^{(j_1)},Q^{(j_3)}]$ and $[Q^{(j_2)},Q^{(j_4)}]$ intersect,\\ because the map $P^{(j)}\mapsto Q^{(j)}$ is linear. Thus there exist $\lambda, \mu \in (0,1)$ such that 
$$
\lambda Q^{(j_1)} + (1-\lambda)Q^{(j_3)} = \mu Q^{(j_2)} + (1-\mu)Q^{(j_4)}.
$$
Therefore, by Proposition~\ref{lem:psizero} and the linearity of $\psi$, we obtain
\begin{align*}
\varphi\left(W_{\{j_1,j_3\}}\right) = \{0\} &\Leftrightarrow 
\psi(Q^{(j_1)}) = \psi(Q^{(j_3)}) = 0\\ &\Leftrightarrow \lambda\psi(Q^{(j_1)}) + (1-\lambda)\psi(Q^{(j_3)}) = 0\\
&\Leftrightarrow \psi\left(\lambda Q^{(j_1)} + (1-\lambda)(Q^{(j_3)}) \right) = 0\\
&\Leftrightarrow \psi\left(\mu Q^{(j_2)} + (1-\mu)(Q^{(j_4)}) \right) = 0\\ &\Leftrightarrow \psi(Q^{(j_2)})=\psi(Q^{(j_4)})=0,
\end{align*}
which, again by Proposition~\ref{lem:psizero}, is equivalent to
\begin{align*}
\varphi\left(W_{\{j_2,j_4\}}\right) = \{0\}.
\end{align*}
\end{proof}

\section*{Acknowledgements}
\noindent
The authors want to thank Krzysztof Bogdan, Kinga Sztonyk and Piotr Zwiernik for their insightful comments. We also thank the anonymous referees for constructive suggestions which helped to improve the paper.

\bibliographystyle{abbrv}
\bibliography{mle}
\end{document}